\documentclass[12pt]{article}

\usepackage{amsmath,amssymb}
\usepackage[mathscr]{eucal}

\newcommand{\ol}{\overline}

\newcommand{\WW}{\mathscr{W}}

\newcommand{\RR}{\mathbb{R}}
\newcommand{\CC}{\mathbb{C}}

\newcommand{\OO}{\mathcal{O}}

\begin{document}
\begin{center}
{\Large Spectral properties and self-adjoint extensions of
	the third power of the radial Laplace operator}\\
\vspace{0.3cm}
T.~A.~Bolokhov\\
\vspace{0.3cm}
{\it St.Petersburg Department of V.\,A.\,Steklov Mathematical Institute\\
         Russian Academy of Sciences\\
         27 Fontanka, St.\,Petersburg, Russia 191023}
\end{center}

\begin{abstract}
    We consider self-adjoint extensions of differential operators
    of the type
$ (-\frac{d^2}{dr^2} + \frac{l(l+1)}{r^2})^3 $
    on the real semi-axis for 
$ l=1,2 $ with
    two kinds of boundary conditions: first that nullify
    the value of a function and its first derivative and second
    that nullify the 4th
($ l=1 $) or the 3rd
($ l=2 $) derivative. 
    We calculate the expressions for the correponding resolvents 
    and derive spectral decompositions. These types of boundary conditions
    are interesting from the physical point of view, especially the
    second ones, which give an example of emergence of long-range action
    in exchange for a singularity at the origin.
\end{abstract}


\section{Introduction}
    The problem of calculation of square root of the sixth-order
    differential operators of the type
\begin{equation}
\label{O3}
    T_{l}^{3} = \bigl(-\frac{d^{2}}{dr^{2}} +\frac{l(l+1)}{r^{2}}\bigr)^{3} ,
	\quad l = 1,2,\ldots ,
\end{equation}
    defined on the real semi-axis
    appears in the construction of the quantum states
    of a free transverse field in spherical coordinates.
    Of particular interest is the action of these operators in the
    presence of specific boundary conditions.
    More precisely, domain of operator is suggested
    to consist of the functions which are regular with a factor of
$ r^{-2} $
    --- that is functions which vanish at the origin along with their
    first derivatives, or of the functions
    whose fourth
($l=1$) or third
($l=2$) 
    derivative vanish at the origin.
    We argue that with either choice of the boundary conditions
    symmetric operator
(\ref{O3})
    has non-trivial self-adjoint extensions.
    We construct the expressions for the resolvent kernels 
    and find the corresponding spectral projectors and densities
    of spectral decompositions.

    As concluding remarks, we discuss basic properties of the
    resulting spectral densities and provide some motivation for our
    specific choices of boundary conditions.

    Operators similar to
(\ref{O3})
    but of the fourth order
    were studied in
\cite{Lapl},
    whereby the corresponding formulas of their connection
    to the transverse
    Laplace operator in the spherical coordinates are given.
    We extensively exploit the latter technique for studying
    the boundary conditions of solutions of the corresponding
    differential equations.

\section{Boundary conditions and symmetricity}
    The differential operation
(\ref{O3})
    defines a symmetric operator on the set of
    sufficiently smooth square integrable functions vanishing at zero along with
    their first five derivatives
\begin{equation}
\label{SP0}
    \WW_{0}^{6} = \{u(r): \int (|u(r)|^{2}+|u^{(6)}(r)|^{2})\,dr < \infty ,
	\quad u(0) = u'(0) = \ldots u^{(5)}(0) = 0 \} .
\end{equation}
    Let us try to enlarge this set and find adequate boundary conditions
    such that formal symmetricity of the operator
(\ref{O3})
    stays intact.
    Define
\begin{equation*}
    B = \int_{0}^{\infty} u(r)T_{l}^{3} \bar{v}(r) dr -
	\int_{0}^{\infty} T_{l}^{3} u(r) \bar{v}(r) dr ,
\end{equation*}
    where
$ T_{l} $
    as above
    is a differential operation of the second order
\begin{equation*}
    T_{l} = -\frac{d^{2}}{dr^{2}} + \frac{l(l+1)}{r^{2}} 
\end{equation*}
    (we are not specifying its domain yet), while
$ u(r) $ and
$ v(r) $
    are some smooth functions which satisfy the integrability conditions of
(\ref{SP0}).
    Then formal symmetricity of the operator defined by 
(\ref{O3})
    requires that
$ B = 0 $.
    A straightforward calculation for it the following expresssion
\begin{multline}
    B = \bigl(
	\frac{3L}{r^{4}}(L-8)(u\bar{v}'-u'\bar{v}) +
\label{BB}
    \frac{6L}{r^{3}}(u\bar{v}''-u''\bar{v})
	+ \frac{3L}{r^{2}}(u'\bar{v}''-u''\bar{v}'+u'''\bar{v}-u\bar{v}''')+\\
+ u\bar{v}^{(5)} - u^{(5)}\bar{v} + u''\bar{v}''' - u'''\bar{v}''
	-u'\bar{v}^{(4)} + u^{(4)}\bar{v}'
	\bigr)\bigr|_{r=0} , \quad L=l(l+1).
\end{multline}
    We dropped here the value at infinity due to decrease of functions
$ u(r) $ and
    $ v(r) $
    along with their derivatives.
    The factors that multiply the powers
$ r^{-4} $,
$ r^{-3} $,
$ r^{-2} $
    in
(\ref{BB})
    depend on
$ r $ themselves.
    In order to have
$ B $
    vanishing, one is to expand those factors in powers of
$ r $,
    substitute them into
$ B $
    and then force all coefficients of non-decreasing powers
    to zero.
    Let us define the expansions of
$ u(r) $ and
$ v(r) $
    as
\begin{equation*}
    u(r) = u(0) + \sum_{j\geq 1} u^{(j)}(0) \frac{r^{j}}{j!},\quad
    v(r) = v(0) + \sum_{j\geq 1} v^{(j)}(0) \frac{r^{j}}{j!} ,
\end{equation*}
    and denote
\begin{equation*}
    \omega^{kj} = u^{(k)}(0) \bar{v}^{(j)}(0) - u^{(j)}(0) \bar{v}^{(k)}(0)\,.
\end{equation*}
    Then
\begin{align*}
    u(r)&\bar{v}'(r) - u'(r)\bar{v}(r) = \omega^{01} + \omega^{02}r
	+ (\omega^{03}+\omega^{12}) \frac{r^{2}}{2}
	+ (\omega^{04}+ 2\omega^{13}) \frac{r^{3}}{6}+\\
	&+ (\omega^{05}+ 3\omega^{14}+ 2\omega^{23})\frac{r^{4}}{24} +\ldots, \\
    u(r)&\bar{v}''(r) - u''(r)\bar{v}(r) = \omega^{02}
	+ (\omega^{03}+\omega^{12})r
	+ (\omega^{04}+ 2\omega^{13}) \frac{r^{2}}{2}+\\
	&+ (\omega^{05}+ 3\omega^{14}+ 2\omega^{23})\frac{r^{3}}{6} +\ldots, \\
    u'(r)&\bar{v}''(r) -u''(r)\bar{v}'(r) +u'''(r)\bar{v}(r)-u(r)\bar{v}'''(r) =\\
	&= \omega^{12} - \omega^{03} - \omega^{04}r
	- (\omega^{14}+\omega^{05}) \frac{r^{2}}{2} + \ldots .
\end{align*}
    Substituting these formulas into
(\ref{BB})
    we find the following expansions for 
$ B $:
    for
$ l=1 $
\begin{equation}
\label{FBC1}
    B = -\frac{36\omega^{01}}{r^{4}} -\frac{24\omega^{02}}{r^{3}}
	-\frac{12\omega^{03}}{r^{2}} - \frac{6\omega^{04}}{r}
	-\frac{3}{2}\omega^{05} - \frac{5}{2}\omega^{14} +8\omega^{23} ,
\end{equation}
    and for
$ l=2 $
\begin{equation}
\label{FBC2}
    B = -\frac{36\omega^{01}}{r^{4}} + \frac{36\omega^{12}}{r^{2}}
	+\frac{24\omega^{13}-6\omega^{04}}{r}
	-\frac{7}{2}\omega^{05} + \frac{7}{2}\omega^{14} +16\omega^{23}.
\end{equation}
    From this we conclude that the operator
$ T_{l}^{3} $
    retains its symmetricity for boundary conditions of the type
\begin{equation}
\label{BC1}
    u(0) = u'(0) = 0, \quad \kappa u''(0) = u'''(0),\quad \kappa \in \RR
\end{equation}
    (reality of
$\kappa$
    follows from the conjugation acting on
$ v $),
    \emph{i.\,e.} when acting on the subspace
\begin{align}
\nonumber
    \WW_{\kappa}^{1} = \{u(r):&\int (|u(r)|^{2}+|u^{(6)}(r)|^{2})\,dr < \infty ,
	\\
\label{SP1}
	& \quad u(0) = u'(0)=0,\quad \kappa u''(0) = u'''(0) \} .
\end{align}
    We shall denote the extended operators
$ T_{l1\kappa}^{3} $.
    Another type of boundary conditions which we use below
    demands vanishing of the fourth derivative for 
$ l=1 $
\begin{equation}
\label{BC21}
    u(0) = u^{(4)}(0) = 0, \quad \kappa u''(0) = u'''(0),\quad \kappa \in \RR,
\end{equation}
    and third derivative for
$ l=2 $
\begin{equation}
\label{BC22}
    u'(0) = u'''(0) = u^{(4)}(0) = 0, \quad \kappa^{5} u(0) = u^{(5)}(0),
	\quad \kappa \in \RR .
\end{equation}
    They induce symmetric operators 
$ T_{12\kappa}^{3} $
    on the set
\begin{align}
\nonumber
    \WW_{\kappa}^{21} = \{u(r): &\int (|u(r)|^{2}+|u^{(6)}(r)|^{2})\,dr
	< \infty , \\
\label{SP21}
    &\quad u(0) = u^{(4)}(0)=0, \quad \kappa u''(0) = u'''(0) \} 
\end{align}
    and
$ T_{22\kappa}^{3} $
    on the set
\begin{align}
\nonumber
    \WW_{\kappa}^{22} = \{u(r): &\int (|u(r)|^{2}+|u^{(6)}(r)|^{2})\,dr 
	< \infty ,\\
\label{SP22}
	    &\quad u(0) = u'''(0) = u^{(4)}(0)=0,
	    \quad \kappa^{5} u(0) = u^{(5)}(0) \} ,
\end{align}
    correspondingly.
    Here and in what follows we tag objects (such as
$ T_{l\xi\kappa}^{3} $) related to the boundary
    conditions
(\ref{BC1})
    with an index
$ \xi =1 $,
    and those related to the conditions
(\ref{BC21}) for
$ l=1 $,
    or to
(\ref{BC22}) for
$ l=2 $,
    with an index
$ \xi =2 $.
    
    Note that we have introduced an additional condition
\begin{equation*}
    \frac{d^{4}u}{dr^{4}}\bigr|_{r=0} = 0
\end{equation*}
    in expression
(\ref{BC22})
    in order to set the coefficient
$ 6\omega^{04} $
    at the power 
$ r^{-1} $ in
(\ref{FBC2})
    to zero.
    At first glance this seems like too strong a restriction.
    However, for the class of function we are going to deal with, this
    condition will automatically follow from the
    regularity of the behaviour of 
$ u(r) $
    at zero, which is necessary for the integrability with square
    of the latter.

\section{Properties of the spherical Bessel functions}
    The deficiency indices of the symmetric operator generated by
    the action of
$ T_{l}^{3} $
    on the subspace
$ \WW_{0}^{6} $
    coincide with the corresponding dimensions of the kernels of the 
    adjoint operators
\begin{equation*}
    (T_{l}^{3} \pm i\rho^{6})^{*} , \quad \rho \in \RR^{+} .
\end{equation*}
    In order for us to find these kernels we need a way to
    construct solutions of the corresponding differential equations.
    With this in mind, we provide some basic aspects
    of the theory of spherical Bessel functions (see e.~g.
\cite{Bessel}).

    Let
$ D_{l} $ be a linear differential operation with the action
\begin{equation}
\label{Dw}
    D_{l} w(r) = r^{l+1} \bigl(\frac{1}{r}\frac{d}{dr}\bigr)^{l} \frac{w}{r} .
\end{equation}
    Then, by means of mathematical induction, one can derive the
    following equality
    (a kind of Rayleigh's formula)
\begin{equation}
\label{QD}
    T_{l} D_{l} w(r) = D_{l} T_{0} w(r) = -D_{l} \frac{d^{2}w}{dr^{2}} .
\end{equation}
    We observe that the solutions of equations involving
$ T_{l}^{3} $
    can be constructed from the solutions of equations
    with the 6th order derivative operator
\begin{equation*}
    T_{0}^{3} = -\frac{d^{6}}{dx^{6}} .
\end{equation*}
    Another remarkable property of the operation
$ D_{l} $
    is that the main term of the asymptotic of the result at infinity
    is just
$ l $-th derivative
\begin{equation}
\label{Dinf}
    D_{l} w(r) = \frac{d^{l}w(r)}{dr^{l}} + \OO(\frac{w(r)}{r}),
	\quad r\to \infty.
\end{equation}

    The behaviour of functions of the type
(\ref{Dw})
    near zero can be described starting from
    the action of 
$ D_{l} $
    upon the monomials
\begin{equation*}
    D_{l} r^{k} = (k-2l+1) \ldots (k-1) r^{k-l} .
\end{equation*}
    If some function
$ w(r) $
    can be expanded at 
$ r\to 0 $
    into
\begin{equation*}
    w(r) = w_{0} + w_{1} r + \ldots + w_{7} \frac{r^{7}}{7!}
	+\OO(r^{8}),
\end{equation*}
    then
\begin{align}
\label{A1}
    D_{1} w(r) &= -\frac{w_{0}}{r} + w_{2}\frac{r}{2} + w_{3}\frac{r^{2}}{3}
	+ w_{4} \frac{r^{3}}{8} + w_{5}\frac{r^{4}}{30}
	+ w_{6} \frac{r^{5}}{144}	+ \OO(r^{6}),\\
\label{A2}
    D_{2} w(r) &= \frac{3w_{0}}{r^{2}} -\frac{w_{2}}{2} + w_{4}\frac{r^{2}}{8}
	+w_{5} \frac{r^{3}}{15} + w_{6}\frac{r^{4}}{48}
	+ w_{7}\frac{r^{5}}{210} + \OO(r^{6}), \\
\nonumber
    D_{3} w(r) &= -w_{0}\frac{15}{r^{3}} + w_{2}\frac{3}{2r} -w_{4}\frac{r}{8}
	+ w_{6}\frac{r^{3}}{48} + w_{7}\frac{r^{4}}{105} + \OO(r^{4}) .
\end{align}
    This shows that the operation
$ D_{l} $
    lowers the power of a function by
$ l $,
    while at the same time turning
$ l $
    of its coefficients to zero.
    Thus we find, that some of the conditions
(\ref{BC1}),
(\ref{BC21}),
(\ref{BC22})
    are automatically satisfied for regular functions.
    Precisely, these ones
\begin{equation*}
    u(0) = 0, \quad l =1 ; \quad
    u'(0) = 0, \quad u(r)r = 0, \quad l = 2.
\end{equation*}
    In general case one also has to check that operation
$ D_{l} $,
    when applied to a set of linearly independent functions,
    yields a linearly independent set again.
    For exponents with different periods, which will arise below, this
    is obvious, and we shall not stop on this question in detail.

\section{Solutions of the differential equations and the deficiency indices}
    We now exploit the above formulas in order to describe kernels
    of the adjoint operators, which are square integrable solutions
    of the sixth order equations
\begin{equation}
\label{6oe}
    \bigl(-\frac{d^{2}}{dr^{2}} +\frac{l(l+1)}{r^{2}}\mp i\rho^{6}\bigr)^{3}
	q_{l\pm}(r) = (T_{l}^{3} \mp i\rho^{6}) q_{l\pm}(r) = 0 .
\end{equation}
    Equation
(\ref{QD})
    when applied to the exponents
\begin{equation*}
    w^{k}_{\pm} = \exp\{e^{\mp i\frac{\pi}{12} \mp i\frac{\pi k}{3}}\rho r\},
	\quad k = 0,1,\ldots 5
\end{equation*}
    gives
\begin{multline*}
    T_{l}^{3} D_{l} 
    \exp\{e^{\mp i\frac{\pi}{12} \mp i\frac{\pi k}{3}} \rho r\} =
    D_{l} T_{0}^{3}
    \exp\{e^{\mp i\frac{\pi}{12} \mp i\frac{\pi k}{3}} \rho r\} = \\
=    \pm i \rho^{6} D_{l} 
    \exp\{e^{\mp i\frac{\pi}{12} \mp i\frac{\pi k}{3}} \rho r\} .
\end{multline*}
    Functions
$ D_{l} \exp\{e^{\mp i\frac{\pi}{12} \mp i\frac{\pi k}{3}} \rho r\} $
    grow exponentially at infinity for 
$ k = 0,1,5 $.
    The rest three solutions can be grouped into such linear combinations
\begin{equation}
\label{glc}
    q^{\xi}_{l\pm} = D_{l}\bigl(
	a^{\xi}_{\pm}\exp\{e^{\mp i\frac{3\pi}{4}}\rho r\}
	+ b^{\xi}_{\pm} \exp\{-e^{\mp i\frac{\pi}{12}}\rho r\}
	+ c^{\xi}_{\pm} \exp\{-e^{\mp i\frac{5\pi}{12}}\rho r\} 
    \bigr),
\end{equation}
    which are regular at zero and satisfy the boundary conditions
(\ref{BC1}),
(\ref{BC21}),
(\ref{BC22}) for the corresponding
$ \xi $ and
$ l $.
    In particular, the solutions
\begin{equation*}
    q^{1}_{l\pm} = D_{l}\bigl(
	\exp\{e^{\mp i\frac{3\pi}{4}}\rho r\}
	+ e^{\mp i\frac{2\pi}{3}} \exp\{-e^{\mp i\frac{\pi}{12}}\rho r\}
	+ e^{\pm i\frac{2\pi}{3}} \exp\{-e^{\mp i\frac{5\pi}{12}}\rho r\} 
    \bigr)
\end{equation*}
    satisfy the boundary conditions
(\ref{BC1}),
    while the solutions
\begin{equation*}
    q^{2}_{l\pm} = D_{l}\bigl(
	e^{\mp i\frac{5\pi}{6}}\exp\{e^{\mp i\frac{3\pi}{4}}\rho r\}
	+ \sqrt{3} \exp\{-e^{\mp i\frac{\pi}{12}}\rho r\}
	+ e^{\pm i\frac{5\pi}{6}} \exp\{-e^{\mp i\frac{5\pi}{12}}\rho r\} 
    \bigr)
\end{equation*}
    obey the boundary conditions
(\ref{BC21}) for
$ l=1 $
    and 
(\ref{BC22}) for
$ l=2 $.

    The solutions 
$ q^{1}_{l\pm}(r) $ and
$ q^{2}_{l\pm}(r) $
    for fixed
$ l $ and sign index
    are linearly independent.
    The condition of integrability with square demands the absence
    of negative powers in the expansion in
$ r $
    near zero.
    This, according to 
(\ref{A1}),
(\ref{A2}),
    forces the vanishing of the sum of coefficients
$ a $, $ b $, $ c $
    from the general expression
(\ref{glc}).
    In other words, one non-trivial linear constraint is imposed upon
    the linear subspace formed by the
    three functions vanishing at infinity
\begin{equation*}
    D_{l} \exp\{e^{\mp i\frac{\pi}{12} \mp i\frac{\pi k}{3}} \rho r\} ,
	\quad k=2,3,4.
\end{equation*}
    This means that the linear span of functions
$ q^{1}_{l\pm}(r) $ and
$ q^{2}_{l\pm}(r) $
    for fixed
$ l $ and sign index
    saturates all square-integrable solutions of equations
(\ref{6oe}) ---
    that is, the whole kernel of the corresponding adjoint operator
\begin{equation*}
    \ker(T_{l}^{3} \pm i\rho^{6})^{*} = \{\alpha q^{1}_{l\pm}(r)
	+ \beta q^{2}_{l\pm}(r) , \quad \alpha,\beta \in \CC \}.
\end{equation*}
    Therefore, we find that the deficiency indices of the operators
$ T_{l}^{3} $ at
$ l=1,2 $
    are
$ (2,\,2) $.

    In the case of
$ l=3 $
    the symmetric operator
$ T_{l}^{3} $
    has one-dimensional subspaces of the solutions of equations
(\ref{6oe}),
    that is, its dificiency indices are $(1,1)$.
    The domains of self-adjoint solutions of this operator, however,
    do not allow one to impose extra conditions of the type
(\ref{BC1}),
(\ref{BC21}),
(\ref{BC22}),
    which would present physical interest, and so we shall not
    dwell on this case here. 
    For values of the angular momentum
$ l $
    greater than 3 self-adjoint extensions for
$ T_{l}^{3} $
    are absent altogether.

    Let us now come back to the vanishing of the fourth derivative in the
    boundary condition
(\ref{BC22}).
    All functions under consideration below consist of three or four
    exponents taken with coefficients the sum of which is zero,
    with the operation
$ D_{l} $
    applied on top, \emph{i.~e.}
    they are of the form
\begin{equation}
\label{expexp}
    D_{l} \bigl( \sum_{k} a_{k} \exp\{e_{k}\chi r\}\bigr) ,
	\quad \sum_{k} a_{k} = 0, \quad e_{k}^{6} = 1 .
\end{equation}
    It is obvious that the expansion in
$ r $
    of such a sum of exponents under
$ D_{l} $
    misses all powers multiple of 6.
    After the operation
$ D_{l} $ is applied, that expansion will miss all powers equal to
$ 6-l $
    modulo 6, including, in particular, the fifth power for
$ l=1 $ and
    the fourth power for
$ l=2 $.
    This property remains valid for a sum or integral of such functions
    taken with regular weights.
    So in our further discussion we shall not demand the vanishing
    of the fourth derivative in the boundary conditions
(\ref{BC22}).

\section{Extensions of symmetric operator and essential self-adjointness}
    Since the basis
$ q^{1}_{l\pm}(r) $, 
$ q^{2}_{l\pm}(r) $
    has a non-trivial irrational normalization with heavy coefficients,
    construction of a general self-adjoint extension for
$ T_{l}^{3} $
    via Cayley transform presents quite a computational problem.
    For that reason we shall restrict ourselves
    with an argument that the symmmetric operators
$ T_{l\xi\kappa}^{3} $
    defined above are essentially self-adjoint
    (see the books
\cite{RS}, \cite{Richt} for the theory of self-adjoint extensions).
    It follows from the properties of an adjoint
    operator and from the definition of an extension that
\begin{equation*}
    (T_{l\xi\kappa}^{3} \pm i\rho^{6})^{*} \subset
    (T_{l}^{3} \pm i\rho^{6})^{*} ,
\end{equation*}
    and, in particular, that
\begin{equation*}
    \ker(T_{l\xi\kappa}^{3} \pm i\rho^{6})^{*} \subset
    \ker(T_{l}^{3} \pm i\rho^{6})^{*} .
\end{equation*}
    That is, in order to prove the essential self-adjointness of
$ T_{l\xi\kappa}^{3} $ 
    it suffices to show that neither of the functions
$ q^{1}_{l\pm}(r) $, 
$ q^{2}_{l\pm}(r) $
    is in the kernel
\begin{equation}
\label{kerT}
    \ker(T_{l\xi\kappa}^{3} \pm i\rho^{6})^{*} .
\end{equation}
    For the case 
$ l=1 $
    this can be achieved by presenting quite an arbitrary sample function
$ p(r) $
    from the domain of
$ T_{1\xi\kappa}^{3} $, \emph{i.\,e.} from
(\ref{SP1}),
(\ref{SP21}),
    with only second and third derivatives non-vanishing at zero
\begin{equation*}
    p(0) = p'(0) = p^{(4)}(0) = p^{(5)}(0) = 0 ,
	\quad p'''(0) = \kappa p''(0) \neq 0.
\end{equation*}
    Integration by parts allows one to cast in the scalar product
    the differential operation onto
$ q_{l\pm}^{1,2}(r) $,
    while in the process the boundary terms do not cancel due
    to the presence of a non-zero imaginary part of
$ \pm i\rho^{6} $
    (here one has to use the expansion
(\ref{A1})
    for
$ q(r) $):
\begin{align*}
    \bigl((T_{1\xi\kappa}^{3}&\pm i\rho^{6})^{*}q_{1\pm}^{j}, p\bigr)
	= \int_{0}^{\infty} q_{1\pm}^{j}(r)\,(T_{1}^{3}\mp i\rho^{6})
	    \bar{p}(r) dr = \\
    &= \int_{0}^{\infty} (T_{1}^{3}\mp i\rho^{6})q_{1\pm}^{j}(r)
	    \,\bar{p}(r) dr + 8\bigl(q_{1\pm}''(0)\bar{p}'''(0)
	    - q_{1\pm}'''(0)\bar{p}''(0)\bigr) = \\
    &= 8\bigl(\frac{4}{3}e^{\mp i\frac{\pi}{12}}\rho^{3}\kappa
	+ \frac{3}{4}\rho^{4}\bigr)\bar{p}''(0) \neq 0 , \quad j=1 ,\\
    &= 8\bigl(\frac{4}{3}e^{\pm i\frac{7\pi}{12}}\rho^{3}\kappa
	+ \frac{3}{4}\sqrt{3}e^{\mp i\frac{2\pi}{3}}\rho^{4}\bigr)\bar{p}''(0)
	  \neq 0 , \quad j=2 .
\end{align*}
    Here
$ \rho^{6} $
    is assumed to be real, although the above inequalities hold for any
$ |\arg \rho^{6}|<\pi/2 $,
    \emph{i.\,e.} for any
$ i\rho^{6} $ in the upper half-plane.
    Similar results are also valid for the
$ l=2 $ case.
    This yields that functions
$ q_{l\pm}^{1,2}(r) $
    do not belong to the kernels of the adjoint operators
(\ref{kerT}),
    which implies that the latter kernels are empty, and
    hence, operators
$ T_{l\xi\kappa}^{3} $
    are essentially self-adjoint.

\section{Resolvent kernel}
    The resolvent kernel is a universal tool for constructing
    spectral projectors
    and functional calculus of operators
(see, \emph{e.\,g.} \cite{Resolvent}).
    
    We shall look for the resolvent kernel
$ R(r,s;z) $
    of the operator
$ T_{l\xi\kappa}^{3} $,
    in form of a function which obeys the differential equation
\begin{equation}
\label{DE}
    \bigl((-\frac{d^{2}}{dr^{2}} +\frac{l(l+1)}{r^{2}})^{3} - z^{6}\bigr)
	R(r,s;z) = \delta(r-s) , \quad 0<\arg z <\frac{\pi}{3} ,
\end{equation}
    as well as the boundary conditions
(\ref{BC1}),
(\ref{BC21}) or
(\ref{BC22}),
    while being symmetric in the arguments
$ r $ and
$ s $,
    and exponentially vanishing at infinity in these arguments.
    The functions
\begin{equation*}
    D_{l} \exp\{i e^{i\pi k/3}zr\} , \quad k=0,1,\ldots 5
\end{equation*}
    satisfy the homogeneous equation
\begin{equation}
\label{heq}
    (T_{l}^{3} - z^{6})\,
	D_{l} \exp\{i e^{i\pi k/3}zr\} = 0 .
\end{equation}
    In this section we assume that 
$ T_{l} $ and
$ T_{l}^{3} $
    are just differential operations without
    specific domain.
    Let us denote the three vanishing solutions of the above equation as
\begin{equation*}
    g_{k}(z,r) = D_{l} \exp\{ie^{i\pi k/3}zr\} , \quad k=0,1,2.
\end{equation*}
    We relate to each solution
$ g_{k}(z,r) $
    a growing solution
\begin{equation*}
    d_{k}(z,r) = D_{l} \exp\{-ie^{i\pi k/3}zr\} , \quad k=0,1,2,
\end{equation*}
    in such a way that each pair
$ g_{k} $ and
$ d_{k} $
    satisfies one and the same second-order equation
\begin{equation}
\label{SOC}
    (T_{l} - e^{2\pi ik/3} z^{2})\,g_{k}(z,r) = 0 , \quad
    (T_{l} - e^{2\pi ik/3} z^{2})\,d_{k}(z,r) = 0 .
\end{equation}

    Let us now construct three solutions
$ h_{k}(z,r) $
    of homogeneous equation
(\ref{heq}),
    which would satisfy the corresponding boundary conditions
(\ref{BC1}),
(\ref{BC21}) or
(\ref{BC22}),
    in such a manner that the main terms of their asymptotics at infinity
    would match 
$ d_{k}(r) $
\begin{align}
\label{h0}
    h_{0}(r) &= d_{0}(r) + \alpha_{0} g_{0}(r) + \beta_{0} g_{1}(r)
	+ \gamma_{0} g_{2}(r) ,\\
\label{h1}
    h_{1}(r) &= d_{1}(r) + \alpha_{1} g_{1}(r) + \beta_{1} g_{2}(r)
	+ \gamma_{1} g_{0}(r) ,\\
\label{h2}
    h_{2}(r) &= d_{2}(r) + \alpha_{2} g_{2}(r) + \beta_{2} g_{0}(r)
	+ \gamma_{2} g_{1}(r) .
\end{align}
    For each set of three coefficients
$ \alpha_{k} $,
$ \beta_{k} $,
$ \gamma_{k} $
    (with fixed
$ \xi  $ and $ l $),
    there are three linear boundary conditions.
    Therefore, in a generic case, these coefficients exist and are unique.
    Direct calculations verify this and give, for the boundary conditions
$ \xi = 1 $ and
$ l=1 $
\begin{align*}
    \alpha_{0} &= (2e^{i\frac{5\pi}{6}}\kappa -3z)/p_{1} , &
	\beta_{0} &= 2e^{\frac{7\pi}{6}} \kappa/p_{1} , &
	\gamma_{0} &= 2e^{-i\frac{\pi}{6}} \kappa/p_{1} , \\
    \alpha_{1} &= -3z/p_{1} , &
	\beta_{1} &= 2e^{i\frac{5\pi}{6}} \kappa/p_{1} , &
	\gamma_{1} &= -2i \kappa/p_{1} , \\
    \alpha_{2} &= -(2e^{i\frac{5\pi}{6}}\kappa + 3z)/p_{1} , &
	\beta_{2} &= 2i \kappa/p_{1} , &
	\gamma_{2} &= 2e^{i\frac{7\pi}{6}} \kappa/p_{1} , 
\end{align*}
    where
\begin{equation*}
    p_{1}(z) = 3z + 2e^{i\frac{\pi}{6}} \kappa ,
\end{equation*}
    and for the angular momentum
$ l=2 $:
\begin{align*}
    \alpha_{0} &= (3e^{i\frac{7\pi}{6}}\kappa -3e^{i\frac{\pi}{3}} z)/p_{2} , &
	\beta_{0} &= 2e^{i\frac{2\pi}{6}} z/p_{2} , &
	\gamma_{0} &= -2z/p_{2} , \\
    \alpha_{1} &= 3e^{i\frac{7\pi}{6}}\kappa/p_{2} , &
	\beta_{1} &= 2e^{i\frac{2\pi}{3}} z/p_{2} , &
	\gamma_{1} &= -2e^{i\frac{\pi}{3}} z/p_{2} , \\
    \alpha_{2} &= (3e^{i\frac{7\pi}{6}}\kappa +2e^{i\frac{\pi}{3}}z)/p_{2} , &
	\beta_{2} &= -2z/p_{2} , &
	\gamma_{2} &= -2e^{i\frac{2\pi}{3}} z/p_{2} , 
\end{align*}
    where
\begin{equation*}
    p_{2}(z) = 2z + 3e^{i\frac{\pi}{6}} \kappa .
\end{equation*}
    In the case of the boundary conditions
$ \xi = 2 $,
$ l=1 $
\begin{align*}
    \alpha_{0} &= (z -2e^{-i\frac{\pi}{6}}\kappa)/p_{1} , &
	\beta_{0} &= -2(e^{i\frac{\pi}{3}}z +e^{i\frac{\pi}{6}}\kappa)/p_{1} , &
	\gamma_{0} &= 2(e^{i\frac{2\pi}{3}}z
	    + e^{i\frac{5\pi}{6}} \kappa)/p_{1} , \\
    \alpha_{1} &= -3z/p_{1} , &
	\beta_{1} &= 2(e^{i\frac{\pi}{3}} z 
	    +e^{-i\frac{\pi}{6}} \kappa)/p_{1} , &
	\gamma_{1} &= -2(e^{i\frac{2\pi}{3}}z + i\kappa)/p_{1} , \\
    \alpha_{2} &= (z + 2i\kappa)/p_{1} , &
	\beta_{2} &= -2(e^{i\frac{\pi}{3}} z + i\kappa)/p_{1} , &
	\gamma_{2} &= 2(e^{i\frac{2\pi}{3}}z -e^{i\frac{\pi}{6}}\kappa)/p_{1} , 
\end{align*}
    where
\begin{equation*}
    p_{1}(z) = z + 2e^{i\frac{\pi}{6}} \kappa ,
\end{equation*}
    and for the angular momentum
$ l=2 $
\begin{align*}
    \alpha_{0} &= (2e^{i\frac{\pi}{6}}\kappa^{5}-z^{5})/p_{2} , &
	\beta_{0} &= 2(z^{5} -e^{i\frac{\pi}{6}} \kappa^{5})/p_{2} , &
	\gamma_{0} &= -2(z^{5} + e^{i\frac{5\pi}{6}} \kappa^{5})/p_{2} , \\
    \alpha_{1} &= 3z^{5}/p_{2} , &
	\beta_{1} &= -2(z^{5} +i\kappa^{5})/p_{2} , &
	\gamma_{1} &= 2(e^{i\frac{\pi}{6}} \kappa^{5} -z^{5})/p_{2} , \\
    \alpha_{2} &= -(z^{5} + 2i\kappa^{5})/p_{2} , &
	\beta_{2} &= -2(z^{5} + e^{i\frac{5\pi}{6}}\kappa^{5})/p_{2} , &
	\gamma_{2} &= 2(z^{5} + i\kappa^{5})/p_{2} , 
\end{align*}
    where
\begin{equation*}
    p_{2}(z) = z^{5} + 2e^{i\frac{5\pi}{6}} \kappa^{5} .
\end{equation*}

    Functions
$ h_{k} $ and
$ g_{k} $
    introduced this way can be used to construct the resolvent
\begin{multline}
\label{ResK}
    R(r,s;z) = \frac{1}{3z^{4}W_{0}}\bigl(
	h_{0}(r)g_{0}(s)\theta(s-r) + h_{0}(s)g_{0}(r) \theta(r-s)\bigr) +\\
    + \frac{e^{i\frac{2\pi}{3}}}{3z^{4}W_{1}}\bigl(
	h_{1}(r)g_{1}(s)\theta(s-r) + h_{1}(s)g_{1}(r) \theta(r-s)\bigr) +\\
    + \frac{e^{i\frac{4\pi}{3}}}{3z^{4}W_{2}}\bigl(
	h_{2}(r)g_{2}(s)\theta(s-r) + h_{2}(s)g_{2}(r) \theta(r-s)\bigr) ,
\end{multline}
    where
$ W_{k} $ are the Wronskians
\begin{equation*}
    W_{k}(z) = d_{k}' g_{k} - d_{k} g_{k}'.
\end{equation*}
    Related equations
(\ref{SOC})
    ensure that
$ W_{k} $
    do not depend on
$ r $
    and can be calculated \emph{e.\,g.} at infinity using expression
(\ref{Dinf})
\begin{align}
\label{W1}
    l&=1: & W_{0} &= -2iz^{3},& W_{1} &= 2iz^{3},& W_{2} &= -2iz^{3}, \\
\label{W2}
    l&=2: & W_{0} &= -2iz^{5},& W_{1} &= 2ie^{i\frac{2\pi}{3}} z^{5},&
	W_{2} &= 2ie^{i\frac{\pi}{3}} z^{5} .
\end{align}
    Symmetry of
$ R(r,s;z) $
    in arguments
$ r $ and
$ s $
    is evident, while the boundary conditions at zero and at infinity
    follow from the construction of functions
$ h_{k} $ and
$ g_{k} $.
    In order to test differential equation
(\ref{DE})
    one can substitute the expressions for
$ h_{k} $
    from
(\ref{h0})--(\ref{h2})
    and then split
$ R $
    into four parts
\begin{equation*}
    R = R_{0} + R_{1} + R_{2}+ R_{g} .
\end{equation*}
    The first three parts contain every term that has growing solution
$ d_{k} $
\begin{equation}
\label{Rk}
    R_{k} = \frac{e^{i\frac{2\pi k}{3}}}{3z^{4}W_{k}}\bigl(
	d_{k}(r)g_{k}(s)\theta(s-r)
	    + d_{k}(s)g_{k}(r) \theta(r-s)\bigr),
\end{equation}
    while the last one includes only vanishing terms
$ g_{k} $
\begin{align*}
    &R_{g} = \frac{1}{3z^{4}} \Bigl(
	\sum_{k} \frac{e^{i\frac{2\pi k}{3}}}{W_{k}}
	    \alpha_{k} g_{k}(r) g_{k}(s) +\\
    &+\frac{1}{W_{0}}\bigl(
	(\beta_{0}g_{1}(r)+\gamma_{0}g_{2}(r)) g_{0}(s) \theta(s-r)
	+ (\beta_{0}g_{1}(s)+\gamma_{0}g_{2}(s)) g_{0}(r) \theta(r-s)
	\bigr) \\
    &+\frac{e^{i\frac{2\pi}{3}}}{W_{1}}\bigl(
	(\beta_{1}g_{2}(r)+\gamma_{1}g_{0}(r)) g_{1}(s) \theta(s-r)
	+ (\beta_{1}g_{2}(s)+\gamma_{1}g_{0}(s)) g_{1}(r) \theta(r-s)
	\bigr) \\
    &+\frac{e^{i\frac{4\pi}{3}}}{W_{2}}\bigl(
	(\beta_{2}g_{0}(r)+\gamma_{2}g_{1}(r)) g_{2}(s) \theta(s-r)
	+ (\beta_{2}g_{0}(s)+\gamma_{2}g_{1}(s)) g_{2}(r) \theta(r-s)
	\bigr) \Bigr) .
\end{align*}
    Here we have re-grouped the first line using the relation
\begin{equation}
\label{Step}
    \theta(r-s) + \theta(s-r) = 1 
\end{equation}
    for the Heaviside step function.
    Performing a direct substitution, we determine that
    the following relations hold for all
$ \xi $ and
$ l $, and independently of 
$ z $
\begin{equation*}
    \frac{\beta_{0}(z)}{W_{0}(z)}
	= \frac{e^{i\frac{2\pi}{3}}\gamma_{1}(z)}{W_{1}(z)} ,	\quad
    \frac{\gamma_{0}(z)}{W_{0}(Z)}
	= \frac{e^{i\frac{4\pi}{3}}\beta_{2}(z)}{W_{2}()} , \quad
    \frac{e^{i\frac{2\pi}{3}}\beta_{1}(z)}{W_{1}(z)}
	= \frac{e^{i\frac{4\pi}{3}}\gamma_{2}(z)}{W_{2}(z)} .
\end{equation*}
    Collecting
$ \theta $-functions
    with the same coefficients we find the following expression for 
$ R_{g} $
\begin{align}
\label{Rg}
    R_{g} &= \frac{1}{3z^{4}}
	\sum_{k} \frac{e^{i\frac{2\pi k}{3}}}{W_{k}}
	    \alpha_{k} g_{k}(r) g_{k}(s)
    + \frac{\beta_{0}}{3z^{4}W_{0}}
	\bigl(g_{1}(r)g_{0}(s) + g_{0}(r)g_{1}(s)\bigr) +\\
\nonumber
    &+ \frac{e^{i\frac{2\pi}{3}}\beta_{1}}{3z^{4}W_{1}}
	\bigl(g_{2}(r)g_{1}(s) + g_{1}(r)g_{2}(s)\bigr)
    + \frac{e^{i\frac{4\pi}{3}}\beta_{2}}{3z^{4}W_{2}}
	\bigl(g_{0}(r)g_{2}(s) + g_{2}(r)g_{0}(s)\bigr) .
\end{align}
    This is a smooth and symmetric function in 
$ r $ and
$ s $,
    and it obviously satisfies the homogeneous equation
\begin{equation*}
    (T_{l}^{3} - z^{6})	R_{g}(r,s;z) = 0 .
\end{equation*}

    Functions
$ R_{k} $ 
    have the standard form of resolvents of second-order operators
    constructed from the solutions of equations
(\ref{SOC})
    (with certain boundary conditions, the form of which is not important
    in this case).
    For this reason they satisfy the second-order differential equations
\begin{equation*}
    (T_{l} - e^{i\frac{2\pi k}{3}} z^{2}) R_{k}(r,s;z)
	=\frac{e^{i\frac{2\pi k}{3}}}{3z^{4}} \delta(r-s) .
\end{equation*}
    Splitting the third-power polynomial
$ T_{l}^{3} - z^{6} $
    into \emph{commuting} factors
\begin{equation*}
    (T_{l}^{3}-z^{6}) =
    (T_{l}-z^{2}) (T_{l}-e^{i\frac{2\pi}{3}} z^{2})
	(T_{l}-e^{i\frac{4\pi}{3}} z^{2}) ,
\end{equation*}
    we can use the latter three equations to write
\begin{align*}
    (T_{l}^{3}&-z^{6})(R_{0}+R_{1}+R_{2}) =
	(T_{l}-e^{i\frac{4\pi}{3}}z^{2})(T_{l}-e^{i\frac{2\pi}{3}}z^{2})
	    (T_{l}-z^{2}) R_{0} +\\
	&+ (T_{l}-z^{2})(T_{l}-e^{i\frac{4\pi}{3}}z^{2})
	    (T_{l}-e^{i\frac{2\pi}{3}}z^{2}) R_{1} +\\
	&+ (T_{l}-e^{i\frac{2\pi}{3}}z^{2})(T_{l}-z^{2})
	    (T_{l}-e^{i\frac{4\pi}{3}}z^{2}) R_{2} =\\
    =&\, \frac{1}{3z^{4}}\bigl(
	(T_{l}-e^{i\frac{4\pi}{3}}z^{2})(T_{l}-e^{i\frac{2\pi}{3}}z^{2})
	+e^{i\frac{2\pi}{3}}(T_{l}-z^{2})(T_{l}-e^{i\frac{4\pi}{3}}z^{2}) +\\
	& +e^{i\frac{4\pi}{3}}(T_{l}-e^{i\frac{2\pi}{3}}z^{2})(T_{l}-z^{2})
	\bigr)\delta(r-s)
    = \delta(r-s) .
\end{align*}
    It is this equation that shows that our function
$ R(r,s;z) $
    satisfies 
(\ref{DE}) and the necessary boundary conditions,
    and, therefore, represents the resolvent kernel
    of the self-adjoint operator
$ T_{l\xi\kappa}^{3} $.

\section{Discrete Spectrum}
    The discrete spectrum projector of the operator
$ T_{l\xi\kappa}^{3} $
    is the sum of the residues (in variable
$ z^{6} $) at the poles of the resolvent
(\ref{ResK}).
    However the only possible pole of the latter corresponds to the zero
$ z_{p} $ of the denominator
$ p(z) $
\begin{align*}
    \xi &= 1,&l&=1,& p(z) &= 3z + 2e^{i\frac{\pi}{6}} \kappa ,
	& z_{p} & = -\frac{2}{3} e^{i\frac{\pi}{6}} \kappa,\\
    \xi &= 1,&l&=2,& p(z) &= 2z + 3e^{i\frac{\pi}{6}} \kappa , 
	& z_{p} & = -\frac{3}{2} e^{i\frac{\pi}{6}}  \kappa ,\\
    \xi &= 2,&l&=1,& p(z) &= z + 2e^{i\frac{\pi}{6}} \kappa ,
	& z_{p} & = -2 e^{i\frac{\pi}{6}}  \kappa ,\\
    \xi &= 2,&l&=2& p(z) &= z^{5} + 2e^{i\frac{5\pi}{6}} \kappa^{5} , 
	& z_{p} & = -2^{1/5} e^{i\frac{\pi}{6}}  \kappa ,
\end{align*}
    which falls into the sector
$ 0 < \arg z < \frac{\pi}{3} $
    only when
$ \kappa <0 $.
    In the latter case the projector sought can be written as
\begin{align*}
    P(r,s) &= R(r,s;z)(z_{p}^{6}-z^{6}) |_{z=z_{p}}
	= 6z_{p}^{5} R(r,s;z) (z_{p}-z) |_{z=z_{p}} =\\
	&= 6z_{p}^{5} R_{g}(r,s;z) (z_{p}-z) |_{z=z_{p}} .
\end{align*}
    We have already dropped the parts
$ R_{k} $
    of the resolvent which contain functions growing at infinity,
    due to the regular behaviour of
$ R_{k} $
    in the sector
$ 0 < \arg z < \frac{\pi}{3} $.
    Then one can substitute coefficients
$ \alpha_{k} $,
$ \beta_{k} $,
$ \gamma_{k} $,
    into
(\ref{Rg})
    and rearrange the terms into a full square, taking
    into account the conjugation relations
\begin{equation*}
    \overline{g_{0}(z_{p},r)} = g_{2}(z_{p},r) ,
    \quad \overline{g_{1}(z_{p},r)} = g_{1}(z_{p},r),
\end{equation*}
    thus transforming the projector into canonical form
\begin{equation*}
    P(r,s) = v^{\kappa}(r) v^{\kappa}(s) .
\end{equation*}
    Here
$ v^{\kappa}(r) $ 
    is a normalized real eigenvector determined by
    the following expressions: for
$ \xi=1 $
\begin{align*}
    v_{1}^{\kappa} &= \sqrt{\frac{-3}{2\kappa}}
	D_{1} \bigl(\exp\{\frac{2}{3}\kappa r\}
	+ \exp\{-\frac{2i\pi}{3}+ \frac{2}{3}e^{-i\frac{\pi}{3}}\kappa r\}
    + \exp\{\frac{2i\pi}{3}+ \frac{2}{3}e^{i\frac{\pi}{3}}\kappa r\} \bigr), \\
    v_{2}^{\kappa} &= \sqrt{\frac{-8}{27\kappa^{3}}}
	D_{2} \bigl(\exp\{\frac{3}{2}\kappa r\}
	+ \exp\{-\frac{2i\pi}{3}+ \frac{3}{2}e^{-i\frac{\pi}{3}}\kappa r\}
    + \exp\{\frac{2i\pi}{3}+ \frac{3}{2}e^{i\frac{\pi}{3}}\kappa r\} \bigr),
\end{align*}
    and for
$ \xi=2 $
\begin{align*}
    v_{1}^{\kappa} &= \sqrt{\frac{-1}{2\kappa}}
	D_{1} \bigl( \sqrt{3}\exp\{2\kappa r\}
	- e^{-i\frac{\pi}{6}}\exp\{2e^{-i\frac{\pi}{3}}\kappa r\} 
	- e^{i\frac{\pi}{6}}\exp\{2e^{i\frac{\pi}{3}}\kappa r\} 
	\bigr) ,\\
    v_{2}^{\kappa} &= \sqrt{\frac{-1}{5\cdot 2^{3/5}\kappa^{3}}}
	D_{2} \bigl( \sqrt{3}\exp\{2^{1/5}\kappa r\}
	- e^{-i\frac{\pi}{6}}\exp\{2^{1/5}e^{-i\frac{\pi}{3}}\kappa r\} -\\
	&\quad - e^{i\frac{\pi}{6}}\exp\{2^{1/5}e^{i\frac{\pi}{3}}\kappa r\} 
	\bigr) .
\end{align*}
    These expressions are valid when 
$ \kappa < 0 $.

\section{Continuous Spectrum}
    It follows from the general expression for resolvent kernel
    of self-adjoint operators (see, \emph{e.\,g.}
\cite{Resolvent}),
    that the density of the spectral measure of the operator
$ T_{l\xi\kappa}^{3} $
    can be calculated as the step of the values of the resolvent at
    the edges of the cut.
    In the case of variable
$ z $:
$ 0 < \arg z < \frac{\pi}{3} $
    this step is the difference
\begin{equation*}
    P_{\lambda} d\lambda = \frac{6\lambda^{5}d\lambda}{2\pi i}
	\bigl(R(r,s;\lambda) - R(r,s;e^{i\frac{\pi}{3}}\lambda)\bigr), \quad
    \lambda \in \RR^{+} ,
\end{equation*}
    which allows one to write a formal equality
\begin{equation*}
    T_{l\xi\kappa}^{3} = \int_{0}^{\infty} \lambda^{6} P_{\lambda} d\lambda
	+ (z_{p}^{6} P)_{\kappa<0} .
\end{equation*}

    The expressions for the resolvent differences are quite bulky,
    involving separate calculations for different
$ \xi $ and 
$ l $
    and we therefore will only give here some general constructions
    and the results.
    To begin with, we need to cancel the
$ \theta $-functions in the
$ R_{k} $ terms
(\ref{Rk}).
    For that we can use the following relations
\begin{align}
\label{Wgd1}
    W_{0}(e^{i\frac{\pi}{3}}z) &= W_{1}(z),
	& W_{1}(e^{i\frac{\pi}{3}}z) &= W_{2}(z),
	& W_{2}(e^{i\frac{\pi}{3}}z)&=-W_{0}(z), \\
    g_{0}(e^{i\frac{\pi}{3}}\lambda,r) &= g_{1}(\lambda,r), &
    g_{1}(e^{i\frac{\pi}{3}}\lambda,r) &= g_{2}(\lambda,r), &
    g_{2}(e^{i\frac{\pi}{3}}\lambda,r) &= d_{0}(\lambda,r), \\
\label{Wgd3}
    d_{0}(e^{i\frac{\pi}{3}}\lambda,r) &= d_{1}(\lambda,r), &
    d_{1}(e^{i\frac{\pi}{3}}\lambda,r) &= d_{2}(\lambda,r), &
    d_{2}(e^{i\frac{\pi}{3}}\lambda,r) &= g_{0}(\lambda,r)
\end{align}
    and write
\begin{align*}
    R_{0}(\lambda) &+ R_{1}(\lambda) + R_{2}(\lambda) -
    R_{0}(e^{i\frac{\pi}{3}}\lambda) - R_{1}(e^{i\frac{\pi}{3}}\lambda)
	- R_{2}(e^{i\frac{\pi}{3}}\lambda) = \\
    =&\, \frac{1}{3\lambda^{4}} \Bigl(
	\frac{1}{W_{0}(\lambda)}\bigl(d_{0}(r)g_{0}(s)\theta(s-r)
	    +d_{0}(s)g_{0}(r)\theta(r-s)\bigr) +\\
	&+\frac{e^{i\frac{2\pi}{3}}}{W_{1}(\lambda)}
	    \bigl(d_{1}(r)g_{1}(s)\theta(s-r)
		+d_{1}(s)g_{1}(r)\theta(r-s)\bigr) +\\
	&+\frac{e^{i\frac{4\pi}{3}}}{W_{2}(\lambda)}
	    \bigl(d_{2}(r)g_{2}(s)\theta(s-r)
		+d_{2}(s)g_{2}(r)\theta(r-s)\bigr) -\\
	&-\frac{e^{i\frac{2\pi}{3}}}{W_{1}(\lambda)}
	    \bigl(d_{1}(r)g_{1}(s)\theta(s-r)
		+d_{1}(s)g_{1}(r)\theta(r-s)\bigr) -\\
	&-\frac{e^{i\frac{4\pi}{3}}}{W_{2}(\lambda)}
	    \bigl(d_{2}(r)g_{2}(s)\theta(s-r)
		+d_{2}(s)g_{2}(r)\theta(r-s)\bigr) +\\
	&+\frac{1}{W_{0}(\lambda)}\bigl(g_{0}(r)d_{0}(s)\theta(s-r)
	    +g_{0}(s)d_{0}(r)\theta(r-s)\bigr)
    \Bigr) .
\end{align*}
    In this form, it is quite easy now to get rid of everything that
    contains the growing functions
$ d_{1}(\lambda,r) $ and
$ d_{2}(\lambda,r) $,
    while in the rest of the terms we gather
$ \theta $-functions into unities
(\ref{Step}), and find
\begin{align*}
    R_{0}(\lambda) &+ R_{1}(\lambda) + R_{2}(\lambda) -
    R_{0}(e^{i\frac{\pi}{3}}\lambda) - R_{1}(e^{i\frac{\pi}{3}}\lambda)
	- R_{2}(e^{i\frac{\pi}{3}}\lambda) = \\
    =&\, \frac{1}{3\lambda^{4}W_{0}(\lambda)} 
	\bigl(d_{0}(r)g_{0}(s)+d_{0}(s)g_{0}(r)\bigr) .
\end{align*}
    The latter expression combined with the difference
$ R_{g}(\lambda) -R_{g}(e^{i\frac{\pi}{3}}\lambda) $
    gives
\begin{align*}
    3&\lambda^{4}\bigl(R(r,s;\lambda)-R(r,s;e^{i\frac{\pi}{3}}\lambda)\bigr)
	= \\
&=\frac{1}{W_{0}}\bigl(d_{0}(r)\ol{d_{0}(s)}+g_{0}(r)\ol{g_{0}(s)}\bigr) 
    + \frac{\alpha_{0}(\lambda)}{W_{0}}g_{0}(r)\ol{d_{0}(s)}
    + \frac{\alpha_{2}(e^{i\frac{\pi}{3}}\lambda)}{W_{0}}d_{0}(r)\ol{g_{0}(s)}
	+\\
&+\frac{e^{i\frac{2\pi}{3}}}{W_{1}(\lambda)}
    \bigl(\beta_{1}(\lambda)-\beta_{0}(e^{i\frac{\pi}{3}}\lambda)\bigr)\bigl(
	g_{2}(r)\ol{g_{2}(s)} + g_{1}(r)\ol{g_{1}(s)}\bigr) + \\
&+\frac{e^{i\frac{2\pi}{3}}}{W_{1}}
    \bigl(\alpha_{1}(\lambda)-\alpha_{0}(e^{i\frac{\pi}{3}}\lambda)\bigr)
	g_{1}(r)\ol{g_{2}(s)}
    +\frac{e^{i\frac{4\pi}{3}}}{W_{2}}
	\bigl(\alpha_{2}(\lambda)-\alpha_{1}(e^{i\frac{\pi}{3}}\lambda)\bigr)
	    g_{2}(r)\ol{g_{1}(s)} -\\
&-\frac{e^{i\frac{2\pi}{3}}\gamma_{0}(e^{i\frac{\pi}{3}}\lambda)}{W_{1}}
    \bigl(d_{0}(r)\ol{g_{2}(s)} + g_{1}(r)\ol{g_{0}(s)}\bigr)
    +\frac{\gamma_{0}(\lambda)}{W_{0}}
	\bigl(g_{2}(r)\ol{d_{0}(s)} + g_{0}(r)\ol{g_{1}(s)}\bigr) -\\
&-\frac{e^{i\frac{4\pi}{3}}\beta_{1}(e^{i\frac{\pi}{3}}\lambda)}{W_{2}}
    \bigl(d_{0}(r)\ol{g_{1}(s)} + g_{2}(r)\ol{g_{0}(s)}\bigr)
    +\frac{\beta_{0}(\lambda)}{W_{0}}
	\bigl(g_{1}(r)\ol{d_{0}(s)} + g_{0}(r)\ol{g_{2}(s)}\bigr) ,
\end{align*}
    where we have also taken into account the conjugation relations
\begin{equation*}
    \ol{g_{0}(\lambda,r)} = d_{0}(\lambda,r) , \quad
    \ol{g_{1}(\lambda,r)} = g_{2}(\lambda,r) .
\end{equation*}
    Using the explicit expressions for the Wronskians
(\ref{W1}),
(\ref{W2})
    and coefficients
$ \alpha_{k} $,
$ \beta_{k} $,
$ \gamma_{k} $,
    the resolvent difference can be represented as a ``square'' of an eigenvector of the continuous spectrum
\begin{equation*}
    P_{\lambda} = \frac{6\lambda^{5}}{2\pi i}\bigl(R(r,s;\lambda)
	-R(r,s;e^{i\frac{\pi}{3}}\lambda)\bigr) 
	= u^{\lambda}_{\kappa}(r) u^{\lambda}_{\kappa}(s) ,
\end{equation*}
    where, for
$ \xi = 1 $
\begin{align*}
    u_{1\kappa}^{\lambda}(r) &= \frac{i}{\sqrt{2\pi}\lambda} D_{1} \bigl(
	\exp\{-i\phi_{1}+i\lambda r\} - \exp\{i\phi_{1}-i\lambda r\} +\\
    &\quad + \frac{2\kappa}{|p_{1}|}
    (\exp\{\frac{i\pi}{6}-e^{-i\frac{\pi}{6}}\lambda r\}
	- \exp\{-\frac{i\pi}{6}-e^{i\frac{\pi}{6}}\lambda r\}) \bigr) \\
    u_{2\kappa}^{\lambda}(r) &= \frac{i}{\sqrt{2\pi}\lambda^{2}} D_{1} \bigl(
	\exp\{\frac{i\pi}{6}-i\phi_{2}+i\lambda r\}
	    - \exp\{-\frac{i\pi}{6}+i\phi_{2}-i\lambda r\} +\\
    &\quad + \frac{2\lambda}{|p_{2}|}
    (\exp\{-\frac{i\pi}{6}-e^{i\frac{\pi}{6}}\lambda r\}
	- \exp\{\frac{i\pi}{6}-e^{-i\frac{\pi}{6}}\lambda r\}) \bigr) 
\end{align*}
    and
\begin{equation*}
    e^{2i\phi_{1}} = \frac{p_{1}(\lambda)}{\bar{p}_{1}(\lambda)}
	= \frac{3\lambda+2e^{i\pi/6}\kappa}{3\lambda+2e^{-i\pi/6}\kappa} , \quad
    e^{2i\phi_{2}} = \frac{p_{2}(\lambda)}{\bar{p}_{2}(\lambda)}
	= \frac{2\lambda+3e^{i\pi/6}\kappa}{2\lambda+3e^{-i\pi/6}\kappa} ,
\end{equation*}
    while for
$ \xi = 2 $
\begin{align*}
    u_{1\kappa}^{\lambda}(r) &= \frac{1}{\sqrt{2\pi}\lambda} D_{1} \bigl(
	\exp\{-i\phi_{1}+i\lambda r\} + \exp\{i\phi_{1}-i\lambda r\} -\\
    &\quad - 2\frac{\kappa+e^{i\frac{\pi}{6}}\lambda}{|p_{1}|}
    \exp\{\frac{i\pi}{6}-e^{-i\frac{\pi}{6}}\lambda r\}
	- 2\frac{\kappa+e^{-i\frac{\pi}{6}}\lambda}{\lambda |p_{1}|}
	\exp\{-\frac{i\pi}{6}-e^{i\frac{\pi}{6}}\lambda r\} \bigr), \\
    u_{2\kappa}^{\lambda}(r) &= \frac{i}{\sqrt{2\pi}\lambda^{2}} D_{2} \bigl(
	\exp\{-i\phi_{2}+i\lambda r\} - \exp\{i\phi_{2}-i\lambda r\} -\\
&\quad - 2\frac{\lambda^{5}-e^{i\frac{\pi}{6}}\kappa^{5}}{|p_{2}|}
    \exp\{-e^{-i\frac{\pi}{6}}\lambda r\}
	- 2\frac{\lambda^{5}-e^{-i\frac{\pi}{6}}\kappa^{5}}{\lambda^{5} |p_{2}|}
	\exp\{-e^{i\frac{\pi}{6}}\lambda r\} \bigr) ,
\end{align*}
    here
\begin{equation*}
    e^{2i\phi_{1}} = \frac{p_{1}(\lambda)}{\bar{p}_{1}(\lambda)}
	= \frac{\lambda+2e^{i\pi/6}\kappa}{\lambda+2e^{-i\pi/6}\kappa} , \quad
    e^{2i\phi_{2}} = \frac{p_{2}(\lambda)}{\bar{p}_{2}(\lambda)}
= \frac{\lambda^{5}-2e^{-i\pi/6}\kappa^{5}}{\lambda^{5}-2e^{i\pi/6}\kappa^{5}} .
\end{equation*}

\section{Conclusion and Discussions}
    As we have already argued above, the real functions
    on the three-dimensional space
    which make up the transverse fields, are related
    to the domains of the operators
$ T_{l\xi\kappa}^{3} $
    via the coefficient of
$ r^{-2} $.
    Therefore the boundary conditions
(\ref{BC1}) for 
$ \xi=1 $
    (the vanishing of the value of the function and its first derivative)
    in terms of
    the original problem is just a requirement of the regularity and
    finite behaviour at the origin.
    
    The boundary conditions for
$ \xi=2 $,
    in terms of the representation of a function in the form
(\ref{expexp}),
    imply the absence of the fifth power in the expansion
    in
$ r $
    of the sum of exponents.
    Then, when acting on such a sum with
$ D_{l} $,
    for 
$ l=1 $
    the fourth power will be absent, while for
$ l=2 $ there will be no third power.
    Rephrased in the terms of general form coefficients
$ a_{k} $
    from expression
(\ref{expexp})
    the latter requirement reads as
\begin{equation*}
    \sum_{k} e_{k}^{5} a_{k} = \sum_{k} e_{k}^{-1} a_{k} = 0.
\end{equation*}
    For a quadratic form of the operator
$ T_{l\xi\kappa}^{3} $
    taken to the power
$ 1/2 $
    this condition may imply the presence in the physical question of
    the transverse functions with slow decrease at infinity.
    In other words, introduction of the dimensional parameter
$ \kappa $
    and singularity in the origin,
    allow for the appearance of a long-range action in the problem.

    The spectral decomposition
\begin{equation}
\label{FA}
    u_{l}^{\lambda}(r) \simeq \frac{\sqrt{2}}{\sqrt{\pi}\lambda^{l}}
	D_{l}(\sin\lambda r)
\end{equation}
    which can be found from that of an essentially self-adjoint operator
\begin{equation*}
    T_{l} = -\frac{d^{2}}{dr^{2}} + \frac{l(l+1)}{r^{2}}, \quad l\geq 1,
\end{equation*}
    by consequent raising of the spectral parameter to the third power,
    corresponds the extensions
$ \kappa = 0 $ for
$ l=1 $ and
$ \kappa = +\infty $
    for
$ l=2 $ ---
    all for the boundary conditions
$ \xi=1 $.
    For non-trivial values of 
$ \kappa $
    the expansion
(\ref{FA})
    is only found in an asymptotic limit at large
$ \lambda $ with
$ l=1 $,
    or at small
$ \lambda $ with
$ l=2 $:
\begin{align*}
    u_{1\kappa}^{\lambda} &\simeq
	\frac{i}{\sqrt{2\pi}\lambda}
    D_{1}(e^{i\lambda r} - e^{-i\lambda r}),\quad \lambda\to\infty\\
    u_{2\kappa}^{\lambda} &\simeq
	\frac{i}{\sqrt{2\pi}\lambda^{2}}
    D_{2}(e^{-i\lambda r} - e^{i\lambda r}),\quad \lambda\to 0 .
\end{align*}
    The other limits, that is, small
$ \lambda $ with 
$ l=1 $,
    and large
$ \lambda $ with
$ l=2 $,
    give mutually similar, but less trivial asymptotics:
\begin{align}
    u_{1\kappa}^{\lambda} &\simeq \frac{i}{\sqrt{2\pi}\lambda}
	D_{1}(e^{-i\frac{\pi}{6}+i\lambda r}
	- e^{i\frac{\pi}{6}-i\lambda r}
	+ e^{i\frac{\pi}{6}-e^{-i\pi/6}\lambda r}
	- e^{-i\frac{\pi}{6}-e^{i\pi/6}\lambda r}),
	\quad \lambda\to 0,\\
\label{ulim2}
    u_{2\kappa}^{\lambda} &\simeq \frac{i}{\sqrt{2\pi}\lambda^{2}}
	D_{2}(e^{-i\frac{\pi}{6}-i\lambda r}
	- e^{i\frac{\pi}{6}+i\lambda r} 
	+ e^{i\frac{\pi}{6}-e^{-i\pi/6}\lambda r}
	- e^{-i\frac{\pi}{6}-e^{i\pi/6}\lambda r}),
	\quad \lambda\to\infty .
\end{align}
    Therefore, for the case of boundary conditions
$ \xi=1 $
    we have a certain duality between the behaviour of
    the spectral density at small
$ \lambda $, $ l=1 $
    and large
$ \lambda $, $ l=2 $,
    and vice versa.

    Different self-adjoint extensions
$ T_{l\xi\kappa}^{3} $
    with boundary conditions
$ \xi=1 $ and
$ \xi=2 $
    have a single intersection (a common extension): for angular momentum
$ l=2 $
    the extensions with parameters
$ \kappa=0 $, 
$ \xi=1 $ match that with
$ \kappa = +\infty $,
$ \xi=2 $
    and they correspond to a spectral density given by expression
(\ref{ulim2}), which now is exact (the discrete spectrum is absent).
    Such an extension has the characteristic properties of
    the absence of a dimensional parameter, regular behaviour 
    in the origin with coefficient
$ r^{-2} $
    and permissibility of a slow decrease at infinity.

    Concluding this work, we can say that we have constructed
    the spectral expansions of self-adjoint extensions
    of a differential operator of
    the sixth order of the type
(\ref{O3}),
    corresponding to the boundary conditions
(\ref{BC1}),
(\ref{BC21}),
(\ref{BC22})
    for the cases of the orbital momentum
$ l=1 $ and
$ l=2 $.
    All the extensions have continuous spectrum
    of multiplicity one, occupying the real positive semi-axis.
    Besides that, for negative values of the parameter
$ \kappa $
    there is a negative eigenvalue corresponding to one-dimensional
    eigenspace.

\section*{Acknowledgements}
    The author is grateful to P.~Bolokhov for the discussion which led
    to the idea of the present study and to the 
    St.~Petersburg department of V.~A.~Steklov Mathematical Institute for
    a favourable working environment.

\end{document}